\documentclass[twoside, 11pt]{article}

\usepackage{mathrsfs,amsfonts,amsmath}
\RequirePackage{ifthen,calc}
\usepackage{theorem}
\usepackage[dvips]{color}

\renewcommand{\baselinestretch}{1.0}

 \catcode`@=11
 \def\@evenhead{\hbox to\textwidth{\footnotesize\rm\thepage \hfill
  {\it Revisions and Responses}}} 

 \def\@oddhead{\hbox to \textwidth{\footnotesize{\it
 } \hfill\thepage}}

 \renewcommand{\section}{\makeatletter
 \renewcommand{\@seccntformat}[1]{{\csname the##1\endcsname.}\hspace{0.45em}}
 \makeatother \@startsection
{section}
{1}
{0pt}
{\baselineskip}
{0.5\baselineskip}
{\normalsize\bfseries\mathversion{bold}}}

\renewcommand{\subsection}{\makeatletter
 \renewcommand{\@seccntformat}[1]{{\csname the##1\endcsname.}\hspace{0.45em}}
 \makeatother \@startsection
{subsection}
{1}
{0pt}
{\baselineskip}
{0.5\baselineskip}
{\normalsize\bfseries\mathversion{bold}}}

\catcode`@=12

\theoremheaderfont{\normalfont\bfseries}
\theorembodyfont{\slshape} {\theorembodyfont{\rmfamily}} {\theorembodyfont{\rmfamily}
}
{\theorembodyfont{\rmfamily} } {\theorembodyfont{\rmfamily}
}
{\theorembodyfont{\rmfamily} } {\theorembodyfont{\rmfamily}
}
{\theorembodyfont{\rmfamily} } {\theorembodyfont{\rmfamily} } {\theorembodyfont{\rmfamily}
}

 \def\beqlb{\begin{eqnarray}}\def\eeqlb{\end{eqnarray}}
 \def\beqnn{\begin{eqnarray*}}\def\eeqnn{\end{eqnarray*}}
 \numberwithin{equation}{section}
\begin{document}
\begin{center}
\renewcommand{\baselinestretch}{1.5}\baselineskip 15pt
{\LARGE\bf  Convergence rate in precise asymptotics for
 Davis law of large numbers}\\
\vspace{1 cm} Lingtao Kong,\\
School of Statistics, Shandong University of Finance and Economics,\\
 Jinan 250014, P. R. China\\\emph{kltgw277519@126.com}
\end{center}

{\bf Abstract:} Let $\{X_n,n\geq 1\}$ be a sequence of i.i.d.
random variables with partial sums $\{S_n,~n\geq 1\}$.
Based on the classical Baum-Katz theorem, a paper by Heyde in 1975
initiated the precise asymptotics for the sum $\sum_{n\geq 1}\mbox{P}(|S_n|\geq\epsilon n)$
as $\epsilon$ goes to zero.
Later, Klesov  determined the convergence rate in
Heyde's theorem. The aim of this paper is to
extend Klesov's result to the precise asymptotics for
Davis law of large numbers, a theorem in Gut and Sp\u{a}taru [2000a].



\vspace{1 cm} \setcounter{page}{1}
\renewcommand{\theequation}{\thesection.\arabic{equation}}
\newcommand{\D}{\displaystyle}
\newcommand{\T}{\textstyle}

\begin{flushleft}
\section{Introduction}
\end{flushleft}
\setcounter{equation}{0}
Let $X,X_1,X_2,\cdots$ be i.i.d. random variables with  partial sums $S_n=\sum_{k=1}^nX_k,~n\geq 1$.
Throughout this paper, we use
 $N$ to denote a standard normal random variable and set $\log n=\ln(x\vee e)$.

Hsu and Robbins [1947] introduced the definition of complete convergence.
Later, Baum and Katz [1965] obtained the following
famous result.

{\bf Theorem BK}~~Let $\{X,X_n,n\geq 1\}$ be a sequence of i.i.d.
random variables with partial sums $\{S_n,n\geq1\}$.
Given $0<p<2$ and $r\geq p$, we have
\begin{equation}
\sum_{n\geq1}n^{\frac{r}{p}-2}
\mbox{P}(|S_n|\geq\epsilon n^{\frac{1}{p}})< \infty,~\mbox{ for any }\epsilon>0,
\end{equation}
if and only if
$\mbox{E}|X|^r<\infty$, and when $r\geq 1$, $\mbox{E}X=0$.
For $r=2, ~p=1$, we refer to Hsu and Robbins [1947] (for the sufficiency)
and Erd\H{o}s [1949, 1950] (for the necessity). For $r=p=1$, Theorem BK
reduces to the theorem of Spitzer [1956].

There have been some extensions of complete convergence in several directions.
The first is to consider the precise rate.
We state the result in Heyde [1975] here.

{\bf Theorem H}~~Let $\{X,X_n,n\geq 1\}$ be a sequence of i.i.d.
random variables with partial sums $\{S_n,n\geq1\}$. We have
\begin{equation}
\lim_{\epsilon\searrow0}\epsilon^2\sum_{n\geq1}
\mbox{P}(|S_n|\geq \epsilon n)=\mbox{E}X^2,
\end{equation}
whenever E$X$=0 and E$X^2<\infty$ (corresponding to the case of $r=2, p=1$ in Theorem BK).
For more results on the precise rate, see Chen [1978], Sp\u{a}taru [1999],
Gut and Sp\u{a}taru [2000a, 2000b], Sp\u{a}taru [2004], Li and Sp\u{a}taru [2012], etc.

The following result of Klesov [1994] studied the rate of convergence
in precise asymptotics of Heyde's result (1.2).

{\bf Theorem K}~~Let $\{X,X_n,n\geq 1\}$ be a sequence of i.i.d.
random variables with partial sums $\{S_n,n\geq1\}$.\\
(1) If $X$ is normal with mean $0$ and variance $\sigma^2>0$.
Then, we have
\begin{equation*}
\lim_{\epsilon\searrow0}
\Big(\sum_{n\geq1}\mbox{P}(|S_n|\geq \epsilon n)-\frac{\sigma^2}{\epsilon^2}\Big)=-\frac{1}{2}.
\end{equation*}
(2) If E$X$=0, E$X^2=\sigma^2>0$, and E$|X|^3<\infty$, then
\begin{equation*}
\lim_{\epsilon\searrow0}
\epsilon^{\frac{3}{2}}\Big(\sum_{n\geq1}
\mbox{P}(|S_n|\geq \epsilon n)-\frac{\sigma^2}{\epsilon^2}\Big)=0.
\end{equation*}

Recently, Gut and Steinebach [2012, 2013] extended Theorem K to Theorems 1, 2 in
Gut and Sp\u{a}taru [2000a],
in which  they got the precise asymptotics
 in  Baum-Katz law of large numbers.
Furthermore, Gut and Sp\u{a}taru [2000a, Theorem 3]
 obtained the precise asymptotics
in Davis law of large numbers as follows.

{\bf Theorem GS}~~Let $\{X,X_n,n\geq 1\}$ be a sequence of i.i.d.
random variables with partial sums $\{S_n,n\geq1\}$.
Suppose that E$X$=0 and E$X^2=\sigma^2<\infty$.
We have, for $0\leq \delta\leq 1$,
\begin{equation}
\lim_{\epsilon\searrow0}\epsilon^{2\delta+2}
\sum_{n\geq 1}\frac{(\log n)^\delta}{n}\mbox{P}(|S_n|\geq \epsilon\sqrt{n\log n})
=\frac{\mbox{E}|N|^{2\delta+2}}{\delta+1}\sigma^{2\delta+2}.
\end{equation}

Motivated by the above results, the purpose of
this paper is to extend Klesov's result to Theorem GS.
We first introduce some notations.

{\bf Definition 1.1}~~Given $0\leq\delta\leq1$, define
\begin{equation}
\gamma_{n,\delta}
:=\sum_{j=1}^n\frac{(\log j)^\delta}{j}-\frac{(\log n)^{\delta+1}}{\delta+1}, ~n\geq1,
\end{equation}
and
\begin{equation}
\gamma_\delta
:=\lim_{n\rightarrow\infty}\gamma_{n,\delta}.
\end{equation}

{\bf Remark 1.1}~~For any $0\leq\delta\leq 1$,
we will prove the convergence of the sequence $\{\gamma_{n,\delta},~n\geq 1\}$ in Section 2.

{\bf Theorem 1.1}~~Let $\{X,X_n,n\geq 1\}$ be a sequence of i.i.d.
random variables with E$X$=0, E$X^2=\sigma^2>0$ and
$S_n=\sum_{k=1}^nX_k,~n\geq 1$.
Given $0\leq \delta\leq 1$, suppose that
E$|X|^{2+\delta}<\infty$ when $0<\delta\leq1$,
and E$\big[|X|^{2}\log(1+|X|)\big]<\infty$ when $\delta=0$.
We have
\begin{equation*}
\lim_{\epsilon\searrow 0}\Big[\sum_{n\geq1}
\frac{(\log n)^{\delta}}{n}
\mbox{P}\Big(|S_n|\geq \epsilon \sqrt{n\log n}\Big)-\frac{\mbox{E}|N|^{2\delta+2}}{\delta+1}
\sigma^{2\delta+2}\epsilon^{-(2\delta+2)}\Big]
=\gamma_\delta-\eta_\delta,
\end{equation*}
where $\gamma_\delta$ is defined in (1.5) and
\begin{equation}
\eta_\delta:=\sum_{n\geq 1}\frac{(\log n)^{\delta}}{n}\mbox{P}(S_n=0).
\end{equation}

{\bf Remark 1.2} We will prove $\eta_\delta<\infty$ in Section 3.

Throughout this paper, we use $C$ to denote positive constants whose values
will be different in different places.
Let $N$ denote a standard normal variable and
$\log n=\ln(x\vee e)$.

Without loss of generality, we assume $\sigma=1 $ in the following two sections.

\section{The normal case}
\setcounter{equation}{0}

We first prove the convergence of the sequence  $\{\gamma_{n,\delta},n\geq 1\}$, where $\gamma_{n,\delta}$ is defined in (1.4).

{\bf Lemma 2.1} Given $0\leq \delta\leq 1$,
the limit of the sequence $\{\gamma_{n,\delta},~n\geq 1\}$ exists.

{\bf Proof:}
By the definition of $\gamma_{n,\delta}$
and the mean value theorem, we have, for $n\geq 1$,
\begin{equation}
\gamma_{n+1,\delta}-\gamma_{n,\delta}=
\frac{\big(\log (n+1)\big)^\delta}{n+1}-\frac{\big(\log(n+1)\big)^{\delta+1}-(\log n)^{\delta+1}}{\delta+1}
=\frac{\big(\log (n+1)\big)^\delta}{n+1}-\frac{(\log \xi)^{\delta}}{\xi},
\end{equation}
for some $\xi\in(n,n+1)$.
Let $g(x)=\frac{(\log x)^\delta}{x}$, then
$g^{\prime}(x)=\frac{(\log x)^{\delta-1}(\delta-\log x)}{x^2}$. Thus, when $x>e^{\delta}$,
$g(x)$ is a decreasing function.
So we have, when $n>e^{\delta}$, $\gamma_{n+1,\delta}-\gamma_{n,\delta}<0$,
which means that $\{\gamma_{n,\delta}, n\geq e^{\delta}\}$ is a decreasing sequence.

By Euler-Maclaurin sum formula (see Cram\'{e}r [1946, p. 124]), we have
\begin{eqnarray*}
\sum_{j=1}^{n}\frac{(\log j)^\delta}{j}&=&\int_{1}^{n}\frac{(\log x)^{\delta}}{x}dx
+\frac{(\log n)^\delta}{2n}-\int_{1}^n P_1(x)d\frac{(\log x)^{\delta}}{x}\\
&=&\frac{(\log n)^{\delta+1}}{\delta+1}
+\frac{(\log n)^\delta}{2n}-\int_{1}^n P_1(x)\frac{(\log x)^{\delta-1}(\delta-\log x)}{x^2}dx,
\end{eqnarray*}
where $P_1(x)=[x]-x+\frac{1}{2}$, and $[x]$ denotes the maximal integer less than or equal to $x$.
Thus,
\begin{eqnarray*}
|\gamma_n|&=&\Big|\frac{(\log n)^{\delta}}{2n}-
\int_{1}^nP_1(x)\frac{(\log x)^{\delta-1}(\delta-\log x)}{x^2}dx
\Big|\\
&\leq&\frac{(\log n)^{\delta}}{2n}+\frac{3}{2}\int_{1}^n\Big|
\frac{(\log x)^{\delta-1}(\delta-\log x)}{x^2}\Big|dx.
\end{eqnarray*}
Since $\int_{1}^{\infty}\frac{(\log x)^{\delta}}{x^2}dx<\infty $, we obtain that
the sequence $\{\gamma_{n,\delta}, n \geq1\}$ has a lower bound.
The monotone bounded theorem shows that
$\lim_{n\rightarrow\infty}\gamma_{n,\delta}$ exists, for any $0\leq\delta\leq1$.

{\bf Proposition 2.1}
Given $0\leq \delta\leq 1$,
we have
\begin{equation*}
\lim_{\epsilon\searrow 0}\Big[\sum_{n\geq1}
\frac{(\log n)^{\delta}}{n}
\mbox{P}\Big(|N|\geq \epsilon \sqrt{\log n}\Big)-\frac{\mbox{E}|N|^{2\delta+2}}{\delta+1}
\epsilon^{-2(\delta+1)}\Big]
=\gamma_\delta,
\end{equation*}
where $N$ denotes a standard random variable.

{\bf Proof:}
Note that
\begin{eqnarray}
&&\sum_{j\geq1}\frac{(\log j)^\delta}{j}\mbox{P}(|N|\geq\epsilon\sqrt{\log j})\nonumber\\
&=&\sqrt{\frac{2}{\pi}}\sum_{j\geq1}\frac{(\log j)^{\delta}}{j}
\int_{\epsilon\sqrt{\log j}}^{\infty}e^{-\frac{y^2}{2}}dy\nonumber\\
&=&\sqrt{\frac{2}{\pi}}\sum_{j\geq1}\frac{(\log j)^{\delta}}{j}
\sum_{n\geq j}\int_{\epsilon\sqrt{\log n}}^{\epsilon\sqrt{\log (n+1)}}e^{-\frac{y^2}{2}}dy\nonumber\\
&=&\sqrt{\frac{2}{\pi}}\sum_{n\geq1}\Big[\sum_{j=1}^{n}\frac{(\log j)^{\delta}}{j}\Big]
\int_{\epsilon\sqrt{\log n}}^{\epsilon\sqrt{\log (n+1)}}e^{-\frac{y^2}{2}}dy\nonumber\\
&=&
\sqrt{\frac{2}{\pi}}\sum_{n\geq1}\Big(\frac{(\log n)^{\delta+1}}{\delta+1}+\gamma_{\delta}\Big)
\int_{\epsilon\sqrt{\log n}}^{\epsilon\sqrt{\log (n+1)}}e^{-\frac{y^2}{2}}dy\nonumber\\
&&~~~~~~~~~~~~~~+\sqrt{\frac{2}{\pi}}\sum_{n\geq1}(\gamma_{n,\delta}-\gamma_{\delta})
\int_{\epsilon\sqrt{\log n}}^{\epsilon\sqrt{\log (n+1)}}e^{-\frac{y^2}{2}}dy\\
&=:&\Delta(\epsilon)+\Gamma(\epsilon),
\end{eqnarray}
where (2.2) follows from the definition of $\gamma_{n,\delta}$.

The mean value theorem for integration shows that there exists $\xi\in(n,n+1)$
such that
\begin{eqnarray}
\int_{\epsilon\sqrt{\log n}}^{\epsilon\sqrt{\log (n+1)}}e^{-\frac{y^2}{2}}dy
&=&\epsilon e^{-\frac{\epsilon^2\log\xi}{2}}(\sqrt{\log (n+1)}-\sqrt{\log n})\nonumber\\
&=&\epsilon e^{-\frac{\epsilon^2\log n}{2}}\big(\frac{n}{\xi}\big)^{\frac{\epsilon^2}{2}}
(\sqrt{\log (n+1)}-\sqrt{\log n}).
\end{eqnarray}
By Taylor expansion, we have
\begin{equation}
\big(\frac{n}{\xi}\big)^{\frac{\epsilon^2}{2}}=1+\frac{\epsilon^2}{2}O(\frac{1}{n})
\end{equation}
and
\begin{equation}
\sqrt{\log (n+1)}-\sqrt{\log n}=\frac{1}{2n\sqrt{\log n}}+O(\frac{1}{n^2}).
\end{equation}
Thus, by (2.4)$\sim$(2.6), we have
\begin{eqnarray}
\Delta(\epsilon)&=&
\sqrt{\frac{2}{\pi}}\sum_{n\geq1}\frac{(\log n)^{\delta+1}}{\delta+1}
\int_{\epsilon\sqrt{\log n}}^{\epsilon\sqrt{\log (n+1)}}e^{-\frac{y^2}{2}}dy
+\gamma_\delta\nonumber\\
&=&\frac{\sqrt{\frac{2}{\pi}}\epsilon}{2(\delta+1)}\sum_{n\geq1}\frac{(\log n)^{\delta+\frac{1}{2}}}{n}
e^{-\frac{\epsilon^2}{2}\log n}
+\gamma_\delta+\epsilon O(1).
\end{eqnarray}
Furthermore, the change of variable $y=\epsilon\sqrt{\log x}$ yields that
\begin{eqnarray}
&&\frac{\sqrt{\frac{2}{\pi}}\epsilon}{2(\delta+1)}
\int_{1}^{\infty}\frac{(\log x)^{\delta+\frac{1}{2}}}{x}
e^{-\frac{\epsilon^2}{2}\log x}dx\nonumber\\
&=&\frac{\sqrt{\frac{2}{\pi}}\epsilon}{2(\delta+1)}
\int_{0}^{\infty}\frac{2y^{2(\delta+1)}}{\epsilon^{2\delta+3}}
e^{-\frac{y^2}{2}}dy\nonumber\\
&=&\frac{1}{\delta+1}\epsilon^{-2(\delta+1)}\mbox{E}|N|^{2(\delta+1)}.
\end{eqnarray}

It remains to estimate $\Gamma(\epsilon)$. By Lemma 2.1, for any
$\rho>0$, we may choose positive integer $n_1$ such that
$|\gamma_{n,\delta}-\gamma_\delta|<\rho$ for $n>n_1$. We split the sum at $n_1$ into
two parts. Thus, we have
\begin{eqnarray}
\limsup_{\epsilon\searrow0}|\Gamma(\epsilon)|&\leq&
\limsup_{\epsilon\searrow0}\mbox{P}\Big(0<|N|\leq\epsilon\sqrt{\log(n_1+1)}\Big)\nonumber\\
&&~~~~~~~~~+\rho\mbox{P}\Big(|N|\geq\epsilon\sqrt{\log (n_1+1)}\Big)\leq2\rho.
\end{eqnarray}
(2.3), (2.7)$\sim$(2.9) and the arbitrary of $\rho$ imply that
we complete the proof of Proposition 2.1.$\Box$

\section{The remainder term}
\setcounter{equation}{0}

{\bf Lemma 3.1 (Heyde [1967])}
Let $\{X, X_n,n\geq 1\}$ be a sequence of i.i.d. random variables with E$X=0$, E$X^2=\sigma^2<\infty$,
and $S_n=\sum_{k=1}^nX_k$. Then, we have
\begin{equation}
\sum_{n\geq1}n^{-1+\frac{\delta}{2}}\sup_x\Big|\mbox{P}(S_n\leq\sigma\sqrt{n} x)-\mbox{P}(N\leq x)\Big|
<\infty,~0\leq\delta<1,
\end{equation}
if and only if E$|X|^{2+\delta}<\infty$, when $~0<\delta<1$,
and E$\big[|X|^{2}\log(1+|X|)\big]<\infty$, when $\delta=0$.

{\bf  Proposition 3.1}~Under the conditions of Theorem 1.1, we have
\begin{equation}
\lim_{\epsilon\searrow 0}\sum_{n\geq 1}
\frac{(\log n)^{\delta}}{n}
\Big[\mbox{P}(|S_n|\geq \epsilon\sqrt{n\log n})-\mbox{P}(|N|\geq \epsilon\sqrt{\log n})\Big]
=-\eta_\delta,
\end{equation}
where $\eta_\delta$ is defined in (1.6).

{\bf Proof:}~
Let
\begin{eqnarray}
\Lambda_n:=\sup_{x}\Big|\mbox{P}(|S_n|\geq \sqrt{n}x)-\mbox{P}(|N|\geq x)\Big|.
\end{eqnarray}
When $0\leq \delta<1$,
it follows from (3.1) that
\begin{eqnarray}
&&\sum_{n\geq 1}\frac{(\log n)^{\delta}}{n}\Lambda_n\nonumber\\
&\leq&C\sum_{n\geq 1}n^{-1+\frac{\delta}{2}}\sup_x\Big|
\mbox{P}(S_n\leq\sigma\sqrt{n} x)-\mbox{P}(N\leq x)\Big|<\infty.
\end{eqnarray}
When $\delta=1$, then E$|X|^3<\infty$.
We make use of the following large deviation estimate (Nagaev [1965]):
\begin{equation}
\Big|\mbox{P}(S_n< \sqrt{n}x)-\mbox{P}(N< x)\big|\leq\frac{C\mbox{E}|X|^3}{\sqrt{n}(1+|x|^3)}.
\end{equation}
It follows from (3.5) that
\begin{eqnarray}
\sum_{n\geq 1}\frac{(\log n)^{\delta}}{n}\Lambda_n
\leq C\sum_{n\geq 1}\frac{(\log n)^{\delta}\mbox{E}|X|^3}{n^{\frac{3}{2}}(1+|x|^3)}
\leq C\sum_{n\geq 1}\frac{(\log n)^{\delta}}{n^{\frac{3}{2}}} <\infty.
\end{eqnarray}

Define
\begin{eqnarray*}
\Lambda_n(\epsilon):=\mbox{P}(|S_n|\geq \epsilon\sqrt{n\log n})-\mbox{P}(|N|\geq \epsilon\sqrt{\log n}).
\end{eqnarray*}
Note that $|\Lambda_n(\epsilon)|\leq \Lambda_n$.
Furthermore,
\begin{eqnarray*}
\Lambda_n(\epsilon)=\mbox{P}(|N|<\epsilon\sqrt{\log n})-\mbox{P}(|S_n|< \epsilon\sqrt{n\log n}).
\end{eqnarray*}
By using the dominated convergence theorem, we have
\begin{eqnarray}
&&\lim_{\epsilon\searrow0}\sum_{n\geq 1}\frac{(\log n)^{\delta}}{n}\Lambda_n(\epsilon)\nonumber\\
&=&\sum_{n\geq 1}\frac{(\log n)^{\delta}}{n}\lim_{\epsilon\searrow0}\Lambda_n(\epsilon)\nonumber\\
&=&-\sum_{n\geq 1}\frac{(\log n)^{\delta}}{n}\mbox{P}(S_n=0)=-\eta_\delta,
\end{eqnarray}
Where (3.7) follows from the continuity of the distribution of $N$.
Thus, We complete the proof of Proposition 3.1.$\Box$

{\bf Remark 3.1} It can easily be verified that $\eta_\delta$ (defined in (1.6)) is finite. Note that
\begin{eqnarray}
\eta_\delta&=&\sum_{n\geq 1}\frac{(\log n)^{\delta}}{n}\mbox{P}(S_n=0)
\leq\sum_{n\geq 1}\frac{(\log n)^{\delta}}{n}\mbox{P}(|S_n|<\frac{1}{\sqrt{n}})\nonumber\\
&\leq&\sum_{n\geq 1}\frac{(\log n)^{\delta}}{n}\Big|\mbox{P}(|S_n|<\frac{1}{\sqrt{n}})-\mbox{P}(|N|<\frac{1}{n})\Big|
+\sum_{n\geq 1}\frac{(\log n)^{\delta}}{n}\mbox{P}(|N|<\frac{1}{n}).
\end{eqnarray}
The first sum in (3.8) is finite by (3.4) and (3.6).
Furthermore,
\begin{equation*}
\mbox{P}(|N|<\frac{1}{n})=\sqrt{\frac{2}{\pi}}\int_{0}^{\frac{1}{n}}e^{-\frac{t^2}{2}}dt\leq C \frac{1}{n},
\end{equation*}
which means that the second sum in (3.8) is less than or
equal to $C\sum_{n\geq1}\frac{(\log n)^{\delta}}{n^2}$.
We complete the proof of the finiteness of $\eta_{\delta}$.

{\bf Proof of Theorem 1.1:} By Propositions 2.1 and 3.1,
we may obtain Theorem 1.1.$\Box$

{\bf Acknowledgment.}  We are grateful to the editors and the anonymous reviewer for their constructive
comments.

\end{document}